\documentstyle[12pt]{amsart}
\begin{document}
\baselineskip=16pt

\title[Deformations of a Calabi-Yau hypersurface]{Infinitesimal
deformations of a Calabi-Yau hypersurface of the moduli space of
stable vector bundles over a curve}

\author{Indranil Biswas}

\address{School of Mathematics, Tata Institute of
Fundamental Research, Homi Bhabha Road, Bombay 400005, India}

\email{indranil@@math.tifr.res.in}

\author{L. Brambila-Paz}

\address{CIMAT, Apdo. Postal 402, C.P. 36240. Guanajuato, Gto,
M\'exico}

\email{lebp@@fractal.cimat.mx}

\thanks{The second author acknowledges the support
of CONACYT grant 28492-E}

\date{}

\begin{abstract}

Let $X$ be a compact connected Riemann surface of genus $g$,
with $g\geq 2$, and ${\cal M}_{\xi}$ a smooth moduli space of
fixed determinant semistable vector bundles of rank $n$, with
$n\geq 2$, over $X$. Take a smooth anticanonical divisor $D$
on ${\cal M}_{\xi}$. So $D$ is a Calabi-Yau variety. We compute
the number of moduli of $D$, namely
$\dim H^1(D,\, T_D)$, to be $3g-4 + \dim H^0({\cal
M}_{\xi}, \, K^{-1}_{{\cal M}_{\xi}})$. Denote by
$\cal N$ the
moduli space of all such pairs $(X',D')$, namely $D'$ is a
smooth anticanonical divisor on a smooth moduli space of
semistable vector bundles over the Riemann surface $X'$. It
turns out that the Kodaira-Spencer map from the tangent space to
$\cal N$, at the point represented by the pair $(X,D)$, to
$H^1(D,\, T_D)$ is an isomorphism. This is proved under the
assumption that if $g =2$, then $n\neq 2,3$, and if $g=3$, then
$n\neq 2$.

\end{abstract}

\maketitle

\section{Introduction}

Let $X$ be a compact connected Riemann surface of genus $g$,
with $g\geq 2$. Let ${\cal M}_{\xi} := {\cal M}(n,\xi)$ denote
the moduli space of stable vector bundles $E$ of rank $n$, with
$n\geq 2$, over $X$, such that the line bundle ${\bigwedge}^n E$
is isomorphic to a fixed holomorphic line bundle $\xi$ over $X$.
The degree $d= \mbox{deg}(\xi)$ and $n$ are assumed to be
coprime. We also assume that if $g=2$, then $n\neq 2,3$, and if
$g=3$, then $n\neq 2$.

The moduli space ${\cal M}_{\xi}$ is a connected smooth
projective variety over $\Bbb C$, and for
fixed $n$, the moduli space ${\cal M}_{\xi}$ is
isomorphic to ${\cal M}_{{\xi}'}$ if ${\xi}'$ is another
holomorphic line bundle with $\mbox{deg}(\xi) = \mbox{deg}
({\xi}')$. We take $\xi$ to be of the form $L^{\otimes d}$,
where $L$ is a holomorphic line bundle over $X$ such that
$L^{\otimes (2g-2)}$ is isomorphic to the canonical line bundle
$K_X$.

The Picard group $\mbox{Pic}({\cal M}_{\xi})$ is
isomorphic to $\Bbb Z$. The
anticanonical line bundle $K^{-1}_{{\cal M}_{\xi}}$ is
isomorphic to ${\Theta}^{\otimes 2}$, where $\Theta$ is the
ample generator of $\mbox{Pic}({\cal M}_{\xi})$, known as the
{\it generalized theta line bundle}.

Let $D$ be a smooth divisor on ${\cal M}_{\xi}$ such that the
holomorphic line bundle ${\cal O}_{{\cal M}_{\xi}}(D)$ over
${\cal M}_{\xi}$ is isomorphic to $K^{-1}_{{\cal M}_{\xi}}$.
Such a divisor is a connected simply connected smooth projective
variety with trivial canonical line bundle. In other words, $D$
is a Calabi-Yau variety.

If we move the triplet $(X,L,D)$, in the space of all triplets
$(X',L',D')$, where $D'$ is a smooth Calabi-Yau hypersurface on a
moduli space of stable vector bundles, of the above type, over
$X'$, then we get deformations of the complex manifold $D$,
simply by associating the complex manifold $D'$ to any triplet
$(X',L' ,D')$. The Kodaira-Spencer infinitesimal deformation map
for this family gives a homomorphism from the tangent space of the
moduli space of triplets $(X,L,D)$, of the above type, into $H^1(D,\,
T_D)$, the space parametrizing the infinitesimal deformations of
the complex manifold $D$. The main result here, [Theorem 3.3],
says

\medskip
{\it The above Kodaira-Spencer infinitesimal deformation map is
an isomorphism.}
\medskip

Consequently, there is an exact sequence
$$
0 \, \longrightarrow \, \mbox{Hom}(l, H^0 ({\cal M}_{\xi},\,
K^{-1}_{{\cal M}_{\xi}})/l) \, \longrightarrow \,
H^1(D,\, T_D) \, \longrightarrow \, H^1(X,\, T_X) \,
\longrightarrow \, 0 \, , \leqno{(1.1)}
$$
where $l \subset H^0 ({\cal M}_{\xi},\, K^{-1}_{{\cal M}_{\xi}})$
is the one dimensional subspace defined by $D$. The above inclusion
map
$$
\mbox{Hom}(l, H^0 ({\cal M}_{\xi},\,  K^{-1}_{{\cal M}_{\xi}})/l)
\, \longrightarrow \, H^1(D,\, T_D)
$$
corresponds to the deformations of $D$ obtained by moving the
hypersurface of the {\it fixed} variety ${\cal M}_{\xi}$, i.e.,
$X$ is kept fixed, and the projection $H^1(D,\, T_D)
\longrightarrow H^1(X,\, T_X)$ in (1.1) is the forgetful
map from the space of infinitesimal
deformations of the triplet $(X,L,D)$
to the space of infinitesimal deformations of $X$. From the
above exact sequence (1.1) it follows immediately that
$$
\dim H^1(D,\, T_D) \, = \, 3g-4 \, + \, \dim H^0({\cal M}_{\xi},
\, K^{-1}_{{\cal M}_{\xi}}) \, .
$$
We note that the dimension of any $H^0({\cal M}_{\xi},
\, {\Theta}^{\otimes k})$, in particular that of
$H^0({\cal M}_{\xi}, \, K^{-1}_{{\cal M}_{\xi}})$, is given by
the Verlinde formula.

Let ${\cal U}_D$ denote the restriction to $X\times D$ of a
Poincar\'e vector bundle over $X\times {\cal M}_{\xi}$. For any
$x\in X$, the vector bundle over $D$, obtained by restricting
${\cal U}_D$ to $x\times D$, is denoted by $({\cal U}_D)_x$. The
following result is used in the proof of Theorem 3.3.

\medskip
{\it For any $x\in X$, the vector bundle $({\cal U}_D)_x$ is stable
with respect to any polarization on $D$. Moreover, the
infinitesimal deformation map
$$
T_xX \, \longrightarrow \, H^1(D,\, \mbox{Ad}(({\cal U}_D))_x)\, ,
$$
for the family ${\cal U}_D$ of vector bundles over $D$ parametrized
by $X$, is an isomorphism.}
\medskip

This result was proved in \cite[Theorem 2.5]{BB} under the
assumption that $n\geq 3$. Here it is extended to the rank two
case [Theorem 2.1].

\section{Restriction of the universal vector bundle}

We continue with the notation of the introduction.

The anticanonical line bundle $K^{-1}_{{\cal M}_{\xi}} :=
{\bigwedge}^{\rm top} T_{{\cal M}_{\xi}}$ is isomorphic to
${\Theta}^{\otimes 2}$ \cite[page 69, Theorem 1]{Ra}, where the
generalized theta line bundle $\Theta$ is the ample generator of
the Picard group $\mbox{Pic}({\cal M}_{\xi})$; the Picard group
is isomorphic to ${\Bbb Z}$.

Let $D \subset {\cal M}_{\xi}$ be a smooth divisor, satisfying
the condition that the line bundle ${\cal O}_{{\cal
M}_{\xi}}(D)$ is isomorphic to $K^{-1}_{{\cal M}_{\xi}}$. Let
$$
{\tau} \, : \hspace{.1in} D \hspace{.1in} \longrightarrow
\hspace{.1in} {\cal M}_{\xi}
$$
denote the inclusion map. Using the Poincar\'e adjunction formula,
we have $K_D \, \cong \, {\tau}^* K_{{\cal M}_{\xi}}\bigotimes
{\tau}^* {\cal O}_{{\cal M}_{\xi}} (D)$. In view of the assumption
${\cal O}_{{\cal M}_{\xi}}(D) \cong K^{-1}_{{\cal M}_{\xi}}$, the
canonical line bundle $K_D$ is trivial. Since the divisor $D$ is
ample, it is connected. Since the moduli space ${\cal M}_{\xi}$
is simply connected, the divisor $D$ is also simply connected.
Therefore, $D$ is a Calabi-Yau variety.

Fix a Poincar\'e vector bundle $\cal U$ over $X\times {\cal
M}_{\xi}$. In other words, for any $m\in {\cal M}_{\xi}$, the
vector bundle over $X$ obtained by restricting $\cal U$ to
$X\times m$ is represented by the point $m$. Let ${\rm
Ad}({\cal U})$ denote the rank $n^2-1$ vector bundle over
$X\times {\cal M}_{\xi}$ defined by the trace zero
endomorphisms of ${\cal U}$. The vector bundle $(\mbox{Id}_X \times
{\tau})^* {\cal U}$ (respectively, $(\mbox{Id}_X \times {\tau})^*
{\rm Ad}({\cal U})$) over $X\times D$ will be denoted by
${\cal U}_D$ (respectively, ${\rm Ad}({\cal U}_D)$).

For any fixed $x\in X$, let ${\cal U}_x$ denote the vector bundle
over ${\cal M}_{\xi}$ obtained by restricting $\cal U$ to
$x\times {\cal M}_{\xi}$. The vector
bundle over $D$ obtained by restricting ${\cal U}_D$ (respectively,
${\rm Ad}({\cal U}_D)$) to $x\times D$ will be denoted by
$({\cal U}_D)_x$ (respectively, ${\rm Ad}({\cal U}_D)_x$).

Since $H^2(D,\, {\Bbb Z}) \, = \, {\Bbb Z}$, the stability of a
vector bundle over $D$ does not depend on the choice of
polarization needed to define the degree of a coherent sheaf
over $D$.

\medskip
\noindent {\bf Theorem\, 2.1.}\, {\it For any point
$x\in X$, the vector bundle $({\cal U}_D)_x$ over $D$ is stable.
Moreover, the infinitesimal deformation map
$$
T_xX \, \longrightarrow \, H^1(D,\, {\rm Ad}({\cal U}_D)_x)
$$
for the family ${\cal U}_D$ of vector bundles over $D$
parametrized by $X$, is an isomorphism.}
\medskip

{\it Proof.}\, If $n\geq 3$ and also $g\geq 3$, then the theorem
has already been proved in \cite[Theorem 2.5]{BB}.

Take a point $x\in X$.
We start, as in the proof of Theorem 2.5 of \cite{BB}, by
considering the exact sequence
$$
0 \, \longrightarrow \, \mbox{Ad}({\cal U})_x\otimes {\cal
O}_{{\cal M}_{\xi}}(-D) \, \longrightarrow \, \mbox{Ad}({\cal
U})_x \, {\buildrel{F}\over\longrightarrow} \,
{\tau}_*\mbox{Ad}({\cal U}_D)_x \, \longrightarrow \, 0 \, ,
$$
over $X\times {\cal M}_{\xi}$, where $F$ denotes the restriction
map. This yields the long exact sequence
$$
H^1({\cal M}_{\xi}, \, \mbox{Ad}({\cal U})_x\otimes {\cal
O}_{{\cal M}_{\xi}}(-D)) \hspace{.1in} \longrightarrow
\hspace{.1in} H^1({\cal M}_{\xi}, \, \mbox{Ad}({\cal U})_x)
$$
$$
\longrightarrow \hspace{.1in} H^1(D,\, \mbox{Ad}({\cal U}_D)_x)
\hspace{.1in} \longrightarrow \hspace{.1in} H^2({\cal M}_{\xi},
\, \mbox{Ad}({\cal U})_x\otimes {\cal O}_{{\cal M}_{\xi}}(-D))
$$
of cohomologies. If we consider $\cal U$ as a family of vector
bundles over ${\cal M}_{\xi}$ parametrized by $X$, then the
infinitesimal deformation map
$$
T_x X \hspace{.1in} \longrightarrow \hspace{.1in} H^1({\cal
M}_{\xi}, \, \mbox{Ad}({\cal U})_x)
$$
is an isomorphism \cite[page 392, Theorem 2]{NR}. In view of the
above long exact sequence, to prove that the infinitesimal
deformation map is surjective it suffices to
establish the following lemma.

\medskip
\noindent {\bf Lemma\, 2.2.}\, {\it If $i=1,2$, then the
following vanishing of cohomology
$$
H^i({\cal M}_{\xi}, \, {\rm Ad}({\cal U}_x)\otimes {\cal
O}_{{\cal M}_{\xi}}(-D)) \hspace{.1in} = \hspace{.1in} 0
$$
is valid.}
\medskip

{\it Proof of Lemma 2.2.}\, This lemma was proved in \cite[Lemma
2.1]{BB} under the assumption that $n\geq 3$. So in the proof we
will assume that $n=2$ and $g\geq 4$.

Let $p$ denote, as in \cite[Section 3]{BB}, the natural
projection of the projective bundle ${\Bbb P}({\cal U}_x)$ over
${\cal M}_{\xi}$ onto ${\cal M}_{\xi}$.  Let $T^{\mbox{rel}}_p$
denote the relative tangent bundle for the projection $p$ from
${\Bbb P}({\cal U}_x)$ to ${\cal M}_{\xi}$.  Since $R^1p_*
T^{\mbox{rel}}_p =0$ and $p_* T^{\mbox{rel}}_p \cong {\rm Ad}({\cal
U}_x)$, for any $i =0,1,2$, the isomorphism
$$
H^i({\cal M}_{\xi}, \, {\rm Ad}({\cal
U}_x)\otimes {\cal O}_{{\cal M}_{\xi}}(-D)) \, = \,
H^i(U,\, p^*K_{{\cal M}_{\xi}}\otimes T^{\mbox{rel}}_{p})
$$
is obtained from the Leray spectral sequence for the map $p$.

If $E$ is a stable vector bundle of rank two
and degree one over $X$, then the vector bundle $E'$ over $X$
obtained by performing an elementary transformation
$$
0 \, \longrightarrow\, E' \, \longrightarrow\, E \,
\longrightarrow\, L_x \, \longrightarrow\, 0 \, ,
$$
where $L_x$ is a one dimensional quotient of the fiber $E_x$, is
semistable. Therefore, we have a morphism, which we will denote
by $q$, from ${\Bbb P}({\cal U}_x)$ to the moduli space ${\cal
M}_{\xi (-x)}$. Here $\xi (-x)$ denotes the line bundle
$\xi\bigotimes {\cal O}_X(-x)$, and ${\cal M}_{\xi (-x)}$ is the
moduli space of semistable vector bundles over $X$ of rank two
and determinant $\xi (-x)$.

Define $U \subset {\Bbb P}({\cal
U}_x)$ to be the inverse image, under that map $q$, of the
stable locus of ${\cal M}_{\xi (-x)}$.

The line bundle $T^{\mbox{rel}}_{p}$ is isomorphic to the
relative canonical bundle $K^{\mbox{rel}}_{q}$ \cite[page
85]{Ty}, \cite{NR}. Therefore, to prove the lemma it suffices
to show that
$$
H^i(U,\, p^*K_{{\cal M}_{\xi}}\otimes K^{\mbox{rel}}_{q})
\, = \, 0 \, , \leqno{(2.3)}
$$
where $i =0,1,2$.

Using the isomorphism of $T^{\mbox{rel}}_{p}$ with
$K^{\mbox{rel}}_{q}$, from
$$
q^*K_{{\cal M}_{\xi (-x)}}\bigotimes K^{\mbox{rel}}_{q}
 \, \cong \, K_U \, \cong \, p^*K_{{\cal M}_{\xi}} \bigotimes
K^{\mbox{rel}}_{p}  \leqno{(2.4)}
$$
we have
$$
p^*K_{{\cal M}_{\xi}}\bigotimes K^{\mbox{rel}}_{q} \, \cong \,
q^*K_{{\cal M}_{\xi (-x)}} \bigotimes \left(K^{\mbox{rel}}_{q}
\right)^{\otimes 3} \, .
$$

Since the restriction of the line bundle $p^*K_{{\cal M}_{
\xi}} \bigotimes K^{\mbox{rel}}_{p}$ to a fiber of the map
$q$ has strictly negative degree, using the above isomorphism,
and the projection formula, we
have
$$
H^i(U,\, p^*K_{{\cal M}_{\xi}}\otimes K^{\mbox{rel}}_{q})
\, = \, H^{i-1}\left({\cal M}_{\xi (-x)}, \, K_{{\cal M}_{\xi
(-x)}}\bigotimes R^1q_* \left(K^{\mbox{rel}}_{q}\right)^{\otimes
3}\right) \, , \leqno{(2.5)}
$$
where $i = 0,1,2$.

The map $q$ is smooth fibration
${\Bbb C}{\Bbb P}^1$ fibration over an
open subset $U'$ of ${\cal M}_{\xi (-x)}$. The assumption that
the genus of $X$ at least four, ensures that the codimension of
the complement of $U'$ is at least four. Therefore, by using the
Hartog type theorem for cohomology, the isomorphism (2.5) is
established.

Setting $i=0$ in (2.5), we conclude that $H^0(U,\, p^*K_{{\cal
M}_{\xi}} \bigotimes K^{\mbox{rel}}_{q}) = 0$.

The following proposition is needed for our next step.

\medskip
\noindent {\bf Proposition\, 2.6.}\, {\it
Let $W$ be a holomorphic vector bundle of rank two over a complex
manifold $Z$, and let $f : {\Bbb P}(V) \longrightarrow Z$ be the
corresponding projective bundle. Then there are
canonical isomorphisms
$$
R^1f_* K^{\otimes 3}_f \, \cong \, S^4(W)\bigotimes
\left({\bigwedge}^2 W^*\right)^{\otimes 2} \, \cong \,
R^0f_* T^{\otimes 2}_f \, ,
$$
where $K_f$ (respectively, $T_f$) is the relative canonical
(respectively, anticanonical) line bundle.}
\medskip

{\it Proof of Proposition 2.6.}\, To construct the
isomorphisms, let $V$ be a complex vector space of
dimension two. Choosing a basis of $V$, we identify the tangent
bundle $T_{{\Bbb P}(V)}$ with ${\cal O}_{{\Bbb P}(V)}(2)$, and
also obtain an identification of the line ${\bigwedge}^2 V^*$
with $\Bbb C$. Since,
$$
H^0({\Bbb P}(V),\, {\cal O}_{{\Bbb P}(V)}(m))\, = \, S^m(V) \, ,
$$
we have an isomorphism of $H^0({\Bbb P}(V),\, T^{\otimes
2}_{{\Bbb P}(V)})$ with $S^4(V)\bigotimes \left({\bigwedge
}^2 V^*\right)^{\otimes 2}$. Now it is a straight forward
computation to check that this isomorphism is $GL(V)$ invariant,
i.e., it does not depend on the choice of a basis of $V$. Therefore,
this pointwise construction of a canonical isomorphism of vector
spaces induces an isomorphism
$$
R^0f_* T^{\otimes 2}_f \, \cong \, S^4(W)\bigotimes \left({\bigwedge
}^2 W^*\right)^{\otimes 2} 
$$
between vector bundles.

To obtain the other isomorphism in the statement
of the proposition, first note that by
the Serre duality we have $H^0({\Bbb P}(V),\, T^{\otimes
2}_{{\Bbb P}(V)}) = H^1({\Bbb P}(V),\, K^{\otimes
3}_{{\Bbb P}(V)})^*$. Now the canonical identification of
$S^4(W)\bigotimes \left({\bigwedge}^2 W^*\right)^{\otimes 2}$
with its dual, namely $S^4(W^*)\bigotimes \left({\bigwedge}^2
W\right)^{\otimes 2}$, gives the other isomorphism. This completes
the proof of Proposition 2.6.$\hfill{\Box}$
\medskip

The isomorphisms in Proposition 2.6
are {\it canonical isomorphisms}, i.e., they are
compatible with the pull back of $W$ using any map $Z'
\longrightarrow Z$, and furthermore, the isomorphisms are
compatible with substituting $W$ by $W\bigotimes L$, where $L$
is a holomorphic line bundle over $Z$.

Combining Proposition 2.6 with (2.5), and using the projection
formula, we get that if $i = 0,1,2$, then
$$
H^i(U,\, p^*K_{{\cal M}_{\xi}}\bigotimes K^{\mbox{rel}}_{q})
\, = \, H^{i-1}\left({\cal M}_{(-x)}, \, K_{{\cal M}_{\xi
(-x)}}\bigotimes q_* \left(T^{\mbox{rel}}_{q}\right)^{\otimes
2}\right) \leqno{(2.7)}
$$
$$
= \, \hspace{.1in} H^{i-1}\left(U, \, q^* K_{{\cal M}_{\xi
(-x)}}\bigotimes \left(T^{\mbox{rel}}_{q}\right)^{\otimes
2}\right) \, .
$$
Indeed, the first isomorphism in (2.7) is a consequence
of (2.5) and Proposition 2.6, and since $R^1p_*
\left(T^{\mbox{rel}}_{q}\right)^{\otimes 2} = 0$, the second
isomorphism in (2.7) is valid. Although there is no universal
vector bundle over $X\times {\cal M}_{\xi (-x)}$, the
properties of the isomorphism $R^1f_* K^{\otimes 3}_f \cong R^0f_*
T^{\otimes 2}_f$ in Proposition 2.6
that were explained earlier, evidently
ensure that the isomorphism in (2.7) is valid. More precisely,
the pointwise construction of the isomorphism between
$R^1q_* \left(K^{\mbox{rel}}_{q}\right)^{\otimes 3}$ and
$q_* \left(T^{\mbox{rel}}_{q}\right)^{\otimes 2}$ gives an
isomorphism of vector bundles.

Using (2.4), and the earlier mentioned fact that
$T^{\mbox{rel}}_{p} \cong K^{\mbox{rel}}_{q}$, we obtain that
$$
q^* K_{{\cal M}_{\xi
(-x)}}\bigotimes \left(T^{\mbox{rel}}_{q}\right)^{\otimes
2} \, \cong \, p^*K_{{\cal M}_{\xi}}\bigotimes
\left(K^{\mbox{rel}}_{p}\right)^{\otimes 3} \, .
$$

Since the restriction of $\left(K^{\mbox{rel}}_{p}
\right)^{\otimes 3}$ to a fiber of $p$ has
strictly negative degree, we have $p_*\left(K^{\mbox{rel}}_{p}
\right)^{\otimes 3} = 0$. Consequently,
the above isomorphism simplifies the terms in (2.7) to
give the following isomorphism
$$
H^{i-1}\left(U, \, q^* K_{{\cal M}_{\xi
(-x)}}\bigotimes \left(T^{\mbox{rel}}_{q}\right)^{\otimes
2}\right)\, = \, H^{i-2}\left({\cal M}_{\xi},\, K_{{\cal
M}_{\xi}} \bigotimes R^1p_*\left(K^{\mbox{rel}}_{p}\right)^{\otimes
3}\right) \leqno{(2.8)}
$$
where $i = 0,1,2$.

Note that we obtain $H^1(U,\, p^*K_{{\cal M}_{\xi}}\bigotimes
K^{\mbox{rel}}_{q} ) = 0$ by setting $i=1$ in (2.8).

In order to complete the proof of the lemma we need
to show that
$$
H^2(U,\, p^*K_{{\cal M}_{\xi}} \bigotimes K^{\mbox{rel}}_{q})
\, = \, 0 \, . \leqno{(2.9)}
$$

To prove the above statement first observe that using (2.8), and
setting $i=2$, we have the following isomorphism
$$
H^2(U,\, p^*K_{{\cal M}_{\xi}} \bigotimes
K^{\mbox{rel}}_{q}) \, = \, H^0\left({\cal M}_{\xi},\,
K_{{\cal M}_{\xi}}\bigotimes R^1p_*\left(K^{\mbox{rel}}_{p}
\right)^{\otimes 3}\right)\, . \leqno{(2.10)}
$$

Now using Proposition 2.6 we have
$$
H^0\left({\cal M}_{\xi},\,
K_{{\cal M}_{\xi}}\otimes R^1p_*\big(K^{\mbox{rel}}_{p}
\big)^{\otimes 3}\right)
\, = \,
H^0\left({\cal M}_{\xi}, \, K_{{\cal M}_{\xi}} \otimes
S^4({\cal U}_x)\otimes \big({\bigwedge}^2
{\cal U}^*_x\big)^{\otimes 2}\right)\, ,
$$
where ${\cal U}_x$, as defined earlier, is the vector bundle
over ${\cal M}_{\xi}$ obtained by restricting the Poincar\'e
bundle $\cal U$ to the subvariety $x\times {\cal U}_{\xi}
\subset X\times {\cal U}_{\xi}$.

The vector bundle ${\cal U}_x$ is known to be stable.
Consequently, the vector bundle
$$
S^4({\cal U}_x)\bigotimes \left({\bigwedge}^2
{\cal U}^*_x\right)^{\otimes 2}
$$
is semistable. Now, since the vector bundle $S^4({\cal U}_x)
\bigotimes \left({\bigwedge}^2 {\cal U}^*_x\right)^{\otimes 2}$
is the dual of itself, its degree is zero. On the other hand,
the degree of $K_{{\cal M}_{\xi}}$ is strictly negative. From
these it follows that the vector bundle
$$
K_{{\cal M}_{\xi}}\bigotimes S^4({\cal U}_x)\bigotimes
\left({\bigwedge}^2 {\cal U}^*_x\right)^{\otimes 2}
$$
does not admit any nonzero section, since it is semistable of
strictly negative degree. In view of (2.10), this establishes
the assertion (2.9). Therefore, the assertion (2.3) is valid.
This completes the proof of the lemma.$\hfill{\Box}$
\medskip

Since we have established, in Lemma 2.2, the rank two analog
of Lemma 2.1 of \cite{BB}, the proof of the stability of the
vector bundle $({\cal U}_D)_x$ for rank at least three, as given
in \cite[Theorem 2.5]{BB}, is also valid for the rank two case
if $g\geq 4$.

We note that \cite[Theorem 2.5]{BB} was proved under the
assumption that $g\geq 3$. However, the proof
remains valid for $g=2$ if the
condition that the rank is at least four is imposed. Under this
condition, the codimension of the subvariety over which the map
$q$ fails to be smooth and proper is sufficiently large in order
to be able to apply the analog Hartog's theorem,
which has been repeatedly used, for the cohomologies in question.

This completes the proof of Theorem 2.1.$\hfill{\Box}$
\medskip

In view of the above Lemma 2.2, all the results established in
Section 2 of \cite{BB} for rank $n \geq 3$ remain valid for rank
two and $g \geq 4$.

\section{Computation of the infinitesimal deformations}

Let $X$ be a compact connected Riemann surface of genus $g$,
with $g\geq 2$. Take a holomorphic line bundle $\xi$ over $X$ of
degree $d$.  Let ${\cal M}_{\xi} := {\cal M}(n,\xi)$ denote the
moduli space of stable vector bundles $E$ of rank $n$ over $X$,
with ${\bigwedge}^n E = \xi$.  For another line bundle ${\xi}'$
of degree $d$, the variety ${\cal M}(n,{\xi}')$ is isomorphic to
${\cal M}(n,\xi)$. Indeed, if $\eta$ is a line bundle over $X$
with ${\eta}^{\otimes n} = {\xi}'\otimes {\xi}^*$, then the map
defined by $E \longmapsto E\otimes\eta$ is an isomorphism from
${\cal M}(n,\xi)$ to ${\cal M}(n,{\xi}')$. Therefore, we can
rigidify (infinitesimally) the choice of $\xi$ by the following
procedure. Fix a line bundle $L$
of degree one over $X$ such that $L^{\otimes (2-2g)}$ is
isomorphic to the tangent bundle $T_X$.  We fix $\xi$ to be
$L^{\otimes d}$.

We will assume that the integers $n$ and $d$ are coprime, and
$n\geq 2$. We will further assume that if $g=2$ then $n\neq
2,3$, and if $g=3$, then $n\neq 2$.

The above numerical assumptions are made in order to ensure that
the assertion in Theorem 2.1 is valid for ${\cal M}_{\xi}$.

Take a smooth divisor $D$ on ${\cal M}_{\xi}$ such that
${\cal O}_{{\cal M}_{\xi}}(D) = K^{-1}_{{\cal M}_{\xi}}$. 
Consider the exact sequence of sheaves
$$
0 \, \longrightarrow \, {\cal O}_{{\cal M}_{\xi}} \,
\longrightarrow \, {\cal O}_{{\cal M}_{\xi}}(D) \,
\longrightarrow \,  {\tau}_*N_D   \, \longrightarrow \, 0
$$
over ${\cal M}_{\xi}$, where $N_D$ is the normal bundle of the
divisor $D$, and ${\tau}$ is the inclusion map of $D$ into ${\cal
M}_{\xi}$. Since
$$
H^1({\cal M}_{\xi}, \, {\cal O}_{{\cal M}_{\xi}}) \, = \, 0 \, ,
$$
using the exact sequence of cohomologies, the space of sections
$H^0(D,\, N_D)$ gets identified with the quotient vector space
$H^0({\cal M}_{\xi},\, {\cal O}_{{\cal M}_{\xi}}(D))/{\Bbb C}$.

Let $\cal S$ denote the space of all divisors $D'$ on ${\cal
M}_{\xi}$ such that $D'$ is homologous to $D$, i.e., they are
represented by the same element in $H^2({\cal M}_{\xi}, \, {\Bbb
Z})$. Therefore, $\cal S$ is identified with ${\Bbb P}H^0({\cal
M}_{\xi},\, K^{-1}_{{\cal M}_{\xi}})$.
The tangent space to $\cal S$, at the point
$[D'] \in {\cal S}$ representing
a divisor $D'$, has the following identification
$$
T_{[D']}{\cal S} \, = \, H^0(D',\, N_{D'}) \, = \, H^0({\cal
M}_{\xi}, \, {\cal O}_{{\cal M}_{\xi}}(D'))/{\Bbb C} \, .
$$

Let $\cal N$ denote the moduli space of triplets of the form
$(X,L,D)$, where $X$, $L$ and $D$ are as above (the line bundle $L$
is a $(2g-2)$-th root of $K_X$). So $\cal N$ is an open subset of
moduli space
of triplets of the form $(X,L,\alpha)$, where $\alpha$ is
a linear subspace of $H^0({\cal M}_{\xi},\, K^{-1}_{{\cal
M}_{\xi}})$ of dimension one. The space $\cal N$ parametrizes a
family of Calabi-Yau varieties, simply by associating the
Calabi-Yau variety $D$ to any triplet $(X,L, D) \in {\cal N}$.

Take a point $\gamma \,:= \, (X,L,D)$ in the moduli space $\cal N$.
Associated to this family is the homomorphism
$$
F \, : \, T_{\gamma} {\cal N} \, \longrightarrow \, H^1(D,\, T_D)
\leqno{(3.1)}
$$
that maps the tangent space $T_{\gamma}{\cal N}$ of $\cal N$ at
$\gamma$ to the space of infinitesimal deformation of the
complex manifold $D$. In other words,
this homomorphism $F$ sends any tangent vector
$v \in T_{\gamma} {\cal N}$ to the corresponding Kodaira-Spencer
infinitesimal deformation class of $D$ for the above
family parametrized by $\cal N$.

The vector space $T_{\gamma} {\cal N}$ fits naturally into the
short exact sequence
$$
0 \, \longrightarrow \, H^0(D,\, N_D) \, \longrightarrow \,
T_{\gamma}{\cal N} \, \longrightarrow \, H^1(X,\, T_X)  \,
\longrightarrow \, 0 \, , \leqno{(3.2)}
$$
where the projection $T_{\gamma} {\cal N} \longrightarrow
H^1(X,\, T_X)$ corresponds to the forgetful map, which sends any
point $(X',L',D') \in {\cal N}$ to the point represented by $X'$ in
the moduli space of Riemann surfaces; the inclusion $H^0(D,\,
N_D) \longrightarrow T_{\gamma} {\cal N}$
in (3.2) corresponds to the
obvious homomorphism $T_{[D]} {\cal S} \longrightarrow T_{\gamma}
{\cal N}$, where $\cal S$, as before, is ${\Bbb P}H^0({\cal
M}_{\xi},\, K^{-1}_{{\cal M}_{\xi}})$, the space of
anticanonical divisors on ${\cal M}_{\xi}$.

\medskip
\noindent {\bf Theorem\, 3.3.}\, {\it The Kodaira-Spencer
infinitesimal deformation map $F$ constructed in (3.1) is an
isomorphism of the tangent space $T_{\gamma} {\cal N}$ with
$H^1(D,\, T_D)$.}
\medskip

{\it Proof.}\, We start by considering the exact sequence
$$
0 \, \longrightarrow \, T_D \,
\longrightarrow \, {\tau}^*T_{{\cal M}_{\xi}} \,
\longrightarrow \, N_D   \, \longrightarrow \, 0
$$
of vector bundles over $D$, where $N_D$ is the normal bundle of
$D$, and ${\tau}$, as before, is the inclusion map of $D$ into ${\cal
M}_{\xi}$. This gives us the exact sequence
$$
H^0(D,\, {\tau}^*T_{{\cal
M}_{\xi}})  \longrightarrow \,H^0(D,\, N_D)  \longrightarrow
H^1(D,\, T_{D}) \longrightarrow H^1(D,\, {\tau}^*T_{{\cal M}_{\xi}})
\longrightarrow H^1(D,\, N_D) \leqno{(3.4)}
$$
of cohomologies.

Since the canonical line bundle $K_D$ is trivial, and $N_D \cong
{\tau}^* K^{-1}_{{\cal M}_{\xi}}$ is ample, the Kodaira vanishing
theorem gives
$$
H^1(D,\, N_D) \, = \, 0 \, . \leqno{(3.5)}
$$
Therefore, the homomorphism $H^1(D,\, T_{D}) \longrightarrow
H^1(D,\, {\tau}^*T_{{\cal M}_{\xi}})$ in (3.4) is surjective.

Our next aim is to show that
$$
H^0(D,\, {\tau}^*T_{{\cal M}_{\xi}}) \hspace{.1in} =
\hspace{.1in} 0 \, , \leqno{(3.6)}
$$
which would be the first step in turning (3.4) into
the short exact sequence (1.1) that we are seeking.

For that purpose, consider the vector bundle ${\cal U}_D$ over
$X\times D$ obtained by restricting a Poincar\'e bundle.  Let
$\phi$ (respectively, $\psi$) denote the projection of $X\times
D$ to $X$ (respectively, $D$). The vector bundle $R^1 {\psi}_*
{\rm Ad}({\cal U}_D)$ over $D$ is naturally isomorphic to
${\tau}^* T_{{\cal M}_{\xi}}$. Also, ${\psi}_* {\rm Ad}({\cal
U}_D) = 0$, as the vector bundle
$({\cal U}_D)_x$ is stable, hence {\it
simple}, for every $x\in X$ [Theorem 2.1]. The vector
bundle ${\rm Ad}({\cal U}_D)$, as in Section 2, is the subbundle
of ${\rm End}({\cal U}_D)$ consisting of trace zero
endomorphisms.  Now, using the Leray spectral sequence for the
projection $\psi$, the isomorphism
$$
H^0(D,\, {\tau}^*T_{{\cal M}_{\xi}}) \, = \, H^1(X\times D,\,
{\rm Ad}({\cal U}_D))
$$
is obtained.

The vector bundle $({\cal U}_D)_x$ over $D$, defined in Section
2, has been proved to be stable in Theorem 2.1. So, we have
$H^0(D,\, {\rm Ad}(({\cal U}_D)_x)) = 0$ for every $x \in X$.
Consequently, the isomorphism
$$
H^1(X\times D,\, {\rm Ad}({\cal U}_D))\, = \, H^0(X,\,
R^1{\phi}_* {\rm Ad}({\cal U}_D))
$$
is obtained.

Now, from the second part of Theorem 2.1 we have a
natural isomorphism
$$
R^1{\phi}_* {\rm Ad}({\cal V}_D) \hspace{.1in} =
\hspace{.1in} T_X
$$
obtained using the Poincar\'e bundle.
Finally, since $H^0(X,\, T_X) =0$, the assertion in (3.6) is an
immediate consequence of the above isomorphism.

Using (3.5) and (3.6), the exact sequence in (3.4) reduces to
$$
0 \, \longrightarrow \, H^0(D,\, N_D) \, \longrightarrow \,
H^1(D,\, T_{D})\, \longrightarrow \, H^1(D,\, {\tau}^*T_{{\cal
M}_{\xi}}) \, \longrightarrow \, 0 \, . \leqno{(3.7)}
$$

The comparison of (3.7) with (3.2) shows that the next step
has to be computation of $H^1(D,\, {\tau}^*T_{{\cal M}_{\xi}})$.

Consider the short exact sequence 
$$
0 \, \longrightarrow \, T_{{\cal M}_{\xi}}\bigotimes {\cal
O}_{{\cal M}_{\xi}}(-D) \, \longrightarrow \, T_{{\cal M}_{\xi}}
\, \longrightarrow \, {\tau}_*{\tau}^* T_{{\cal M}_{\xi}}   \,
\longrightarrow \, 0
$$
of sheaves over ${\cal M}_{\xi}$. We know that $H^2({\cal
M}_{\xi},\, T_{{\cal M}_{\xi}}) = 0$ \cite[page 391, Theorem
1.a]{NR}. Also, we have (3.6). Consequently, the exact sequence
yields the long exact sequence
$$
0 \, \longrightarrow \, H^1({\cal
M}_{\xi},\, T_{{\cal M}_{\xi}}\bigotimes K_{{\cal M}_{\xi}})\,
\longrightarrow \, H^1({\cal M}_{\xi},\, T_{{\cal M}_{\xi}})
\leqno{(3.8)}
$$
$$
\, \longrightarrow \,
H^1(D,\,  {\tau}^*T_{{\cal M}_{\xi}})\, \longrightarrow \, H^2({\cal
M}_{\xi},\, T_{{\cal M}_{\xi}}\bigotimes K_{{\cal M}_{\xi}})
\, \longrightarrow \, 0
$$
of cohomologies; note that $H^i({\cal M}_{\xi},\, {\tau}_*{\tau}^*
T_{{\cal M}_{\xi}}) = H^i(D,\, {\tau}^*T_{{\cal M}_{\xi}})$.

It was proved in \cite{NR} that the Kodaira-Spencer deformation
map for ${\cal M}_{\xi}$, as the Riemann surface $X$ moves
in the moduli space of Riemann surfaces, is
an isomorphism of $H^1({\cal M}_{\xi},\, T_{{\cal M}_{\xi}})$
with $H^1(X,\, T_X)$. Therefore, comparing (3.2) with (3.7), and
using the exact sequence (3.8), we conclude that in order to
complete the proof of the theorem, it suffices to
establish the following statement~: if $i=1,2$, then
$$
H^i({\cal M}_{\xi},\, T_{{\cal M}_{\xi}}\bigotimes K_{{\cal
M}_{\xi}})\, = \, 0 \, . \leqno{(3.9)}
$$
Indeed, (3.9) implies that $H^1(D,\,  {\tau}^*T_{{\cal
M}_{\xi}}) = H^1({\cal M}_{\xi},\, T_{{\cal M}_{\xi}}) =
H^1(X,\, T_X)$.

To prove (3.9), let $\delta$ denote the dimension of the variety
${\cal M}_{\xi}$. The Serre duality gives the following isomorphism
$$
H^i({\cal M}_{\xi},\, T_{{\cal M}_{\xi}}\bigotimes K_{{\cal
M}_{\xi}})\, = \, H^{\delta -i}({\cal M}_{\xi},\,
{\Omega}^1_{{\cal M}_{\xi}})^* \, = \, H^{1,\delta -i}({\cal
M}_{\xi})^* \, . \leqno{(3.10)}
$$
(Here $H^{j,k}({\cal M}_{\xi}) \, := \, H^k({\cal M}_{\xi},\,
{\Omega}^j_{{\cal M}_{\xi}})$.)

To finish the proof of the statement (3.9) we need to use some
properties of the Hodge structure of the cohomology algebra
$H^*({\cal M}_{\xi}, \, {\Bbb C})$, which will be recalled now. 

Fix a Poincar\'e bundle $\cal U$ over $X\times {\cal M}_{\xi}$.
Let $c_k \, := \, c_k({\cal U}) \, \in \, H^{k,k}(X\times {\cal
M}_{\xi})$ denote the $k$-th Chern class of
$\cal U$. For any ${\alpha} \in H^{i,j}(X)$, we have
$$
\lambda (k,\alpha) \, :=\,  \int_X c_k\cup f^*\alpha \, \in \,
H^{k+i-1,k+j-1}({\cal M}_{\xi})\, , \leqno{(3.11)}
$$
where $f$ denotes the obvious projection of $X\times {\cal M}_{\xi}$
onto $X$, and $\int_X$ is the Gysin map for this projection, which
is constructed by integrating differential forms on
$X\times {\cal M}_{\xi}$ along the fibers of the projection $f$.
The collection of all these cohomology classes $\{\lambda
(k,\alpha)\}$, constructed in (3.11), generate the cohomology
algebra $H^*({\cal M}_{\xi},\, {\Bbb C})$ \cite[page 581, Theorem
9.11]{AB}. On the other hand, we know that the following
$$
H^{0,1}({\cal M}_{\xi}) \hspace{.1in} = \hspace{.1in} 0
$$
is valid.

With these properties of $H^*({\cal M}_{\xi},\, {\Bbb C})$ at
our disposal, we are in a position to prove
that the algebra generated by the cohomology classes
$\{\lambda (k,\alpha)\}$ cannot have a nonzero element in
$H^{1,\delta -1}({\cal M}_{\xi})$ or $H^{1,\delta -2}
({\cal M}_{\xi})$, where $\delta = {\dim}_{\Bbb C} {\cal M}_{\xi}$.

To prove the above assertion, suppose that
$$
{\omega} \, = \, {\omega}_1\wedge {\omega}_2\wedge \cdots \wedge
{\omega}_l
$$
is a nonzero element in $H^{1,\delta -1}({\cal
M}_{\xi})\bigoplus H^{1,\delta -2}({\cal M}_{\xi})$, where
${\omega}_j \in \{\lambda (k,\alpha)\}$ for all $j \in [1,l]$.
We will see that ${\omega}_j \in H^{0,1}({\cal M}_{\xi})$ for at least
one $j \in [1,l]$. Since $H^{0,1}({\cal M}_{\xi}) =0$, this would
prove that $\omega = 0$.

First observe that in (3.11), we
have $k+i-1 \geq k-1$ and $k+j-1 \leq k$, as ${\dim}_{\Bbb C} X =1$.
In other words, we have
$$
(k+j-1) - (k+i-1) \hspace{.1in} \leq \hspace{.1in} 1 \, .
\leqno{(3.12)}
$$
Let ${\omega}_i \in H^{a_i,b_i}({\cal M}_{\xi})$, where $i\in
[1,l]$. Then $a_i \leq 1$, and consequently from (3.12)
the inequality $b_i\leq 2$ is obtained. Furthermore,
$a_j \neq 0$ for at most one $j\in [1,l]$. If $a_i =0$, then
$b_i\leq 1$; but the possibility $b_i =1$ is ruled out as 
$H^{0,1}({\cal M}_{\xi})  = 0$. Therefore, all ${\omega}_i$
except one is a scalar. Now, if $a_j=1$, then from (3.12)
we have $b_j\leq 2$. On the other hand, we have $\delta -2 > 2
\geq b_j$. Consequently, we conclude that $\omega = 0$.

Since the cohomology classes $\lambda (k,\alpha)$ are of pure
type, i.e.,
$$
\lambda (k,\alpha) \hspace{.1in} \in \hspace{.1in}
H^{a,b}({\cal M}_{\xi})
$$
for some integers $a$ and $b$, it is easy to see that for any
$i\geq 0$, the cohomology group $H^{i}({\cal
M}_{\xi},\, {\Bbb C})$ is generated, as a complex vector space,
by completely decomposable elements, i.e., elements of the type
$\omega$ considered above. Therefore, we have $H^{1,\delta
-1}({\cal M}_{\xi}) = 0 = H^{1,\delta -2}({\cal M}_{\xi})$.

In view of (3.10), this completes the proof of the statement
(3.9). We already noted that the statement (3.9) completes the
proof of the theorem.$\hfill{\Box}$
\medskip

As a consequence of Theorem 3.3 we get that
$$
\dim H^1(D,\, T_D) \, = \, 3g-4 \, + \, \dim H^0({\cal M}_{\xi},
\, {\Theta}^{\otimes 2}) \, .
$$
The dimension of $H^0({\cal M}_{\xi}, \, {\Theta}^{\otimes 2})$
is given by the Verlinde formula.

\medskip
\noindent {\bf Remark\, 3.13.}\, From
\cite[page 760, Proposition 1]{Be}, coupled
with \cite[page 759, Th\'eor\`eme 1]{Be}, it follows that
$H^0(D,\, T_D) =0$. We note that this is also an immediate
consequence of (3.6).

We have $H^2(D,\, T_D) \,= \, H^{m-1,2}(D) \, = \,
H^{m,3}({\cal M}_{\xi})$,
where $m= \dim D$; the second isomorphism is
obtained from the Lefschetz hyperplane theorem \cite[page 156]{GH}.
The earlier proof that $H^{1,\delta -i}({\cal M}_{\xi}) =0$
for $i=1,2$, easily extends to prove that
$H^{1,\delta -3}({\cal M}_{\xi}) =0$. Therefore, we have
$H^2(D,\, T_D) =0$. However, by a theorem due to
Bogomolov-Kawamata-Tian-Todorov it is already known that the
deformations of a Calabi-Yau variety are unobstructed.


\end{document}